\documentclass[a4paper]{amsart}

\usepackage{amsfonts}
\usepackage{verbatim}
\usepackage{amssymb}
\usepackage{amscd}
\usepackage{amsmath}
\usepackage{latexsym}
\usepackage{framed, color}

\usepackage[dvips]{graphicx}
\usepackage{psfrag}

\makeatletter
\@namedef{subjclassname@2020}{%
  \textup{2020} Mathematics Subject Classification}
\makeatother

\newtheorem{thm}{Theorem}[section]

\newtheorem{lem}[thm]{Lemma}

\theoremstyle{definition}
\newtheorem{dfn}[thm]{Definition}

\theoremstyle{remark}
\newtheorem{rem}[thm]{Remark}

\newtheorem{ex}[thm]{Example}

\makeatletter
 
 \@addtoreset{equation}{section}
\makeatother

\newcommand{\R}{\mathbb{R}}

\renewcommand{\tilde}{\widetilde}

\def\spmapright#1{\smash{%
 \mathop{\hbox to 1.3cm{\rightarrowfill}}
  \limits^{#1}}}
\def\spmapleft#1{\smash{%
 \mathop{\hbox to 1.3cm{\leftarrowfill}}
  \limits^{#1}}}

\title[Reeb spaces of smooth functions on manifolds II]
{Reeb spaces of smooth functions on manifolds II}
\author{Osamu Saeki} 
\address{Institute of Mathematics for Industry,
Kyushu University,
Motooka 744, Nishi-ku, Fukuoka 819-0395, Japan}

\date{\today}


\begin{document}
\begin{abstract}
The Reeb space of a continuous function
is the space of connected components 
of the level sets.
In this paper we characterize those
smooth functions on closed manifolds whose
Reeb spaces have the structure of a finite
graph.
We also give several
explicit examples of smooth functions on closed
manifolds such that they themselves or their Reeb spaces have some
interesting properties.
\end{abstract}

\maketitle

\section{Introduction}\label{section1}

Let us consider a smooth function $f : M \to \R$
on a smooth closed manifold $M$ of dimension $m \geq 2$.
The \emph{Reeb space} $W_f$ of $f$ is the topological space
consisting of the connected components of the level sets.
Thus, we have the natural quotient map $q_f : M \to W_f$, and then
there exists a unique continuous function $\bar{f} : W_f
\to \R$, called the \emph{Reeb function} of $f$, such that
$f = \bar{f} \circ q_f$ (for precise definitions,
see \S\ref{section2}).

In our previous paper \cite{Sa22}, we showed that
the Reeb space of a smooth function on a closed
manifold with only finitely many critical values
always has the structure of a finite (multi-)graph without loops.
Note that the
Reeb space of such a smooth function has the structure
of a graph in such a way that
the vertices correspond exactly to
the connected components of level sets that contain critical points.
Furthermore,
the restriction of the Reeb function to each edge
of the graph is an embedding into $\R$.
In the literature, a Reeb space of a smooth function
is often called a 
\emph{Reeb graph} (or a \emph{Kronrod--Reeb graph}), and
our theorem justifies the terminology (see \cite{R, Sharko}).
If we impose, on the graph structure, the above additional
conditions, then the finiteness of critical values
is clearly a necessary condition.
However, it is not a necessary
condition for the existence of just a finite graph structure
for the Reeb space, as we see in \S\ref{section3} 
of the present paper.

On the other hand, Gelbukh \cite{G0}--\cite{G5}
has extensively studied the topology of Reeb spaces
of smooth functions on manifolds, or even
those of continuous functions on certain topological spaces.
In particular, in \cite{G5}, it is shown that 
the Reeb space of a smooth function on a closed connected manifold
is always a $1$--dimensional Peano continuum 
(compact connected locally connected Hausdorff topological space),
and that it is homotopy equivalent to a finite graph.

In the present paper, we first characterize
those smooth functions on closed manifolds
whose Reeb spaces have the structure of a finite graph.
For this, we use the topological
characterization of an open interval due to Ward \cite{Ward} and
Franklin--Krishnarao \cite{FK1, FK2}
together with Gelbukh's results mentioned above.

We also give explicit examples of
Reeb spaces of smooth functions on closed manifolds
of dimension $m \geq 2$
that have interesting properties with
respect to critical points, critical values, the number
of level set components, and the Reeb spaces.

The paper is organized as follows. In \S\ref{section2}, we review the
basic definitions and certain properties of Reeb spaces.
After preparing some necessary notions, we state and prove the main
theorem of this paper, characterizing those smooth functions whose
Reeb spaces have the structure of a finite graph.
We also see that the theorem can naturally be generalized
to continuous functions on separable Peano continuum;
i.e., we give a topological characterization of continuous
functions on such spaces whose Reeb spaces have the structure of a 
finite graph.

In \S\ref{section3}, we first give an example
of a smooth function on a closed manifold
which has infinitely many critical values but
whose Reeb space has the structure of a graph.
Then, we give an example of a smooth function on
the $2$--sphere $S^2$ which has an infinite number of connected
components of critical point set but which
has only finitely many critical values.
As the third example, we construct a smooth function
$F : S^2 \to \R$ whose level sets all have exactly two connected
components except for $F^{-1}(0)$, which is connected,
such that the Reeb space $W_F$ does not have the
structure of a graph.
Finally, we give an example of a smooth function on the standard
$m$--dimensional sphere $S^m$, $m \geq 2$,
whose Reeb space has the structure of a graph but such that
the associated Reeb function to $\R$
restricted to an edge is not an embedding.

Throughout the paper, 
all smooth manifolds and smooth maps between them are differentiable
of class $C^\infty$ unless otherwise specified. 
A manifold is \emph{closed}
if it is compact and has no boundary.

\section{Reeb space and its graph structure}\label{section2}

Let $f : X \to Y$ be a continuous map
between topological spaces.
For two points $x_0, x_1 \in X$, we write
$x_0 \sim x_1$ if $f(x_0) = f(x_1)$ and $x_0, x_1$
belong to the same connected component
of $f^{-1}(f(x_0)) = f^{-1}(f(x_1))$.
Let $W_f = X/\!\sim$ be the quotient space 
with respect to this equivalence relation.
Let $q_f : X \to W_f$ denote the quotient map.
Then, there exists a unique continuous 
map $\bar{f} : W_f \to Y$ which makes the 
following diagram commutative:
\begin{eqnarray*}
X \!\!\!\! & \spmapright{f} & \!\!\!\! Y \\
& {}_{q_f}\!\!\searrow \quad \qquad \nearrow_{\bar{f}} & \\
& \,W_f. &
\end{eqnarray*}
The space $W_f$ is called the
\emph{Reeb space} of $f$, and
the map $\bar{f} : W_f \to Y$ is called
the \emph{Reeb map} of $f$ (when $Y = \R$,
it is also called the \emph{Reeb function} of $f$). 
The decomposition
$f = \bar{f} \circ q_f$ as in the above
commutative diagram is called the \emph{Stein
factorization} of $f$ \cite{L1}. 

In this paper, a 
\emph{graph} means a finite ``multi-graph''
which may contain multi-edges and/or loops.
When considered as a topological space, it
is a compact $1$--dimensional
CW complex. 

Let $M$ be a smooth closed manifold of
dimension $m \geq 2$ and $f : M \to \R$
a smooth function.
Then, the following has been proved in \cite{Sa22}.

\begin{thm}\label{thm1}
If $f$ has at most finitely many critical values, then
the Reeb space $W_f$ has the structure
of a graph.
Furthermore, the Reeb function $\bar{f} : W_f \to \R$ is an embedding
on each edge.
\end{thm}

Note, however, that the finiteness of critical values
is not a necessary
condition for the graph structure
of the Reeb space (see Example~\ref{ex1} in \S\ref{section3}).
In order to characterize those smooth functions
whose Reeb spaces have the structure of a graph,
let us prepare some notions.

\begin{dfn}\label{def1}
Let $f : X \to \R$ be a continuous function on a 
topological space.
For each $t$, the space $f^{-1}(t)$ is called
a \emph{level set}, and a component of $f^{-1}(t)$
is called a \emph{contour} (or a \emph{level set component}).
A subset $A \subset X$ is \emph{saturated}
if it is the union of contours of $f$:
in other words, if $q_f^{-1}(q_f(A)) = A$.

Now, suppose that $X$ is a smooth manifold and that $f : X \to \R$
is a smooth function. Then,
a contour is \emph{critical}
if it contains a critical point; otherwise, it is said to
be \emph{regular}.
\end{dfn}

\begin{dfn}\label{def2}
Let $f : M \to \R$ be a smooth function
on a closed manifold $M$ of dimension $m \geq 2$.

(1) A non-empty open connected saturated subset $U \subset M$ is 
\emph{cylindrical}
if for every $z \in q_f(U)$,
$U \setminus q_f^{-1}(z)$ has exactly two connected components.

(2) Let $C$ be a contour.
An open connected saturated neighborhood $U$ of $C$ in $M$
is \emph{star-like} if $U \setminus C$ is a finite disjoint
union of cylindrical open connected sets.
Here, $U$ can be equal to $C$ when $C$ is a component of $M$.
\end{dfn}

\begin{figure}[tbh]
\centering
\psfrag{1}{$(1)$}
\psfrag{2}{$(2)$}
\psfrag{R}{$\R$}
\psfrag{f}{$f$}
\psfrag{U}{$U$}
\psfrag{C}{$C$}
\includegraphics[width=\linewidth,height=0.4\textheight,
keepaspectratio]{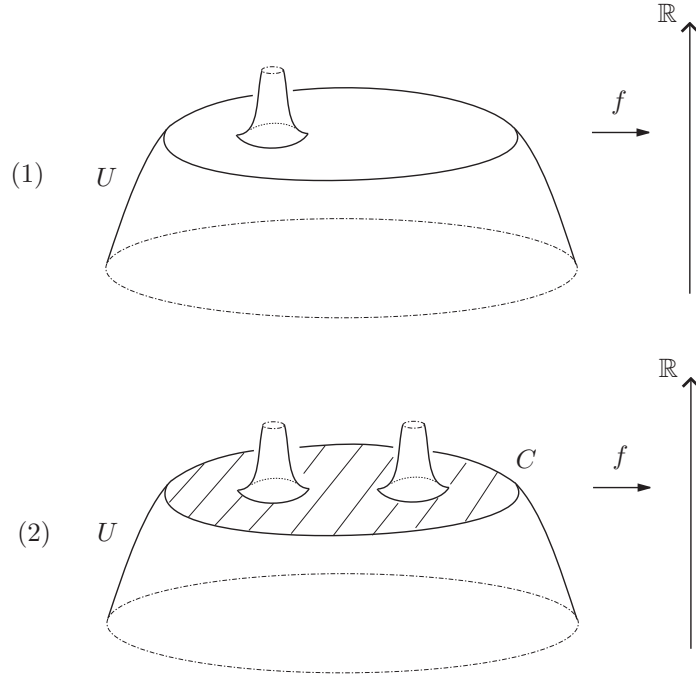}
\caption{(1) Cylindrical open subset \quad (2) Star-like open neighborhood}
\label{fig24}
\end{figure}

For examples of a cylindrical subset and a star-like open
neighborhood, refer to Fig.~\ref{fig24}.
The singular fiber of a Morse function on a nonorientable
closed surface as in \cite[Fig.~2.2 (3)]{Sa04}
has a cylindrical open neighborhood, while the ones
in \cite[Fig.~2.2 (1), (2)]{Sa04} have star-like open neighborhoods.

The following is proved in \cite{G0}.
For the reader's convenience, we give a sketch of proof here.

\begin{lem}\label{lem1}
Let $f : M \to \R$ be a smooth function
on a closed manifold of dimension $m \geq 2$. Then,
every regular contour has a cylindrical
open connected neighborhood.
\end{lem}

\begin{proof}
Let $R$ be a regular contour and set $f(R) = r \in \R$.
As the set of regular points is open in $M$,
by \cite[Proposition~14]{G5}, there exists an 
$\varepsilon > 0$ such that the connected component $U$
of $f^{-1}(I)$ containing $R$ is open, saturated and consists
of regular points for
$I = (r - \varepsilon, r + \varepsilon)$.
Then the function $f|_U : U \to I$ is a proper
submersion.
Therefore, by
an elementary argument using gradient-like vector fields and their
integral curves, we see that $U$ is diffeomorphic
to $R \times I$ and that $f|_U$ corresponds to
the projection to the second factor $R \times I \to I$.
Thus, $U$ is the desired cylindrical open connected 
neighborhood of $R$.
\end{proof}

\begin{lem}\label{lem2}
If $U \subset M$ is a cylindrical open connected set, then
$q_f(U)$ is homeomorphic to an open interval.
\end{lem}

\begin{proof}
Since $U$ is an open subset of $M$, it is separable and
hence so is $q_f(U)$.
The space $q_f(U)$ is Hausdorff, since so is $W_f$ by \cite{G5}.
Furthermore, for every point $z \in q_f(U)$,
$U \setminus q_f^{-1}(z)$ consists of two connected components,
which are open and saturated. Since
$q_f(U) \setminus \{z\}$ is the disjoint union
of their $q_f$--images, which are open, we see that
$q_f(U) \setminus \{z\}$ consists of
exactly two connected components.

Then, by \cite{Ward}, \cite{FK1, FK2} (see also
\cite[\S28]{Willard}), $q_f(U)$ is homeomorphic to
an open interval. (Note that $W_f$ and hence $q_f(U)$
is metrizable, see \cite[Theorem~18]{G5}.)
\end{proof}

\begin{rem}\label{rem0.5}
The above proof shows that if $U$
is a cylindrical open connected subset of $M$, then
for every contour $L$ contained in $U$, $U \setminus
L$ is the disjoint union of two cylindrical open connected sets.
So, $U$ is a star-like neighborhood of $L$.
\end{rem}

Now, the following characterizes those
smooth functions whose Reeb spaces have the
structure of a finite graph.

\begin{thm}\label{thm2}
Let $f : M \to \R$ be a smooth function
on a closed manifold of dimension $m \geq 2$.
Then, the Reeb space $W_f$ has the structure
of a finite graph
if and only if every critical contour
has a star-like open neighborhood.
\end{thm}

\begin{proof}
First suppose that $W_f$ has the structure
of a graph. Let $C$ be a critical contour of $f$. 
Then, $v = q_f(C) \in W_f$
has a small open connected neighborhood $V$
in $W_f$ which is the union
of $v$ and a finite number of
half open arcs $J_i$ each
homeomorphic to $[0, 1)$, $i = 1, 2, \ldots,
\ell$, attached to $v$ at the
end points corresponding to $0$,
for some $\ell \geq 0$.
Then, $U = q_f^{-1}(V)$ is a saturated
open neighborhood of $C$ in $M$.

Suppose $U$ is not connected.
Then, it is the disjoint union of
non-empty open sets $O_0$ and $O_1$.
Since for each $w \in V$, $q_f^{-1}(w)$
is connected, it is entirely contained in
either $O_0$ or $O_1$. This implies that $O_0$ and
$O_1$ are saturated open sets.
Then, we see that $V$ is the disjoint
union of the non-empty open sets $q_f(O_0)$ and $q_f(O_1)$,
which contradicts the connectedness of $V$.
Therefore, $U$ is connected.
By similar arguments, we can show that
$U \setminus C$ is the disjoint union
of $U_i = q_f^{-1}(J_i \setminus \{v\})$, which
are connected and open. Furthermore, for every
point $w \in q_f(U_i) = J_i \setminus \{v\}$,
$U_i \setminus q_f^{-1}(w)$ consists exactly of
two connected components. Therefore,
$U_i$ is cylindrical.
Consequently, $U$ is the desired star-like open neighborhood of $C$.

Conversely, suppose that every critical 
contour has a star-like open neighborhood.
Note that every regular contour
has a cylindrical open connected neighborhood by
Lemma~\ref{lem1}.
Hence, $M$ is covered by such cylindrical open connected
sets and star-like open neighborhoods
of critical contours. As $M$ is compact, we can
choose finitely many of them so that they cover the whole $M$.
(Here, we warn the reader that a star-like open neighborhood
of a critical contour may contain other critical contours.
Therefore, there might be infinitely many critical contours.)

Let $U$ be a member of the finite family of star-like open
neighborhoods as above. We assume that $U$ is a star-like
open neighborhood of a critical contour $C$. 
Suppose $U \setminus C$ has
$\ell$ connected components $U_i$, $i = 1, 2, \ldots, \ell$. 
Let $S_\ell$
be the topological space obtained from $\ell$ copies
of the half open interval $[0, 1)$ by attaching
them to a point, say $p$, at their end points corresponding to $0$.
%Note that $T_\ell \subset S_\ell$, and that $T_\ell$ is compact.
Let us show that $q_f(U)$ is homeomorphic to $S_\ell$
by taking a smaller $U$ if necessary.

Since each component $U_i$ of $U \setminus C$
is cylindrical, by Lemma~\ref{lem2},
$q_f(U_i)$ is homeomorphic to the open interval $(0, 1)$.

Suppose that the closure of a component $U_i$ in $U$
does not intersect with $C$. In this case, if the closure intersects with
some $U_j$ with $j \neq i$, then this is a contradiction,
since $U_i$ and $U_j$ are distinct connected components
of $U \setminus C$. Then, $U$ would be disconnected having a
component containing $C$ and another component $U_i$.
This is a contradiction. Therefore,
the closure of each $U_i$ in $U$ intersects with $C$.

Then, we have an infinite sequence of points in $U_i$ converging
to a point in $C$. The $q_f$--image of the sequence converges
to $q_f(C)$ contained in the closure of $q_f(U_i)$ in $W_f$,
since $q_f$ is continuous
(note that $M$ and $W_f$ are metrizable, see \cite[Theorem~18]{G5}).
If $q_f(C)$ is a point in $q_f(U_i)$, this is a contradiction,
since $U_i \cap C = \emptyset$. Hence, by identifying
$q_f(U_i)$ with the open interval $(0, 1)$ up to
homeomorphism, we see that the above sequence
can be considered to be a sequence in $[0, 1]$, and by
choosing a subsequence if necessary, we may assume that it
converges either to $0$ or to $1$.
If there are two sequences converging to $0$ and to $1$,
respectively, then by taking off the set corresponding to
$q_f^{-1}([1/3, 2/3])$ from $U$ and $U_i$ and by increasing $\ell$
by $1$, we may assume that for every sequence
of points in $U_i$ converging to a point in $C$, a subsequence of
its $q_f$--image converges to $0$ (and never to $1$).

Let $r_i : q_f(U_i) \cup \{q_f(C)\} \to [0, 1)$ be
the map such that $r_i|_{q_f(U_i)}$ is a homeomorphism
onto $(0, 1)$ as above and $r_i(q_f(C)) = 0$.
The above argument shows that there exists a sequence
of $q_f(U_i)$ converging to $q_f(C)$ whose $r_i$--image converges to $0$.
Then, for any sequence in $q_f(U_i)$ converging to
$q_f(C)$, its $r_i$--image converges to $0$,
since $r_i|_{q_f(U_i)}$ is a homeomorphism.
Hence, $r_i$ is a continuous bijection.

Let us show that its inverse $r_i^{-1} : [0, 1) \to 
q_f(U_i) \cup \{q_f(C)\}$ is also continuous.
Take an arbitrary sequence 
of $(0, 1)$ converging to $0$. Then, its $r_i^{-1}$--image
is a sequence of $q_f(U_i) \subset W_f$.
As $W_f$ is compact, there is a subsequence, say $\{w_n\}$,
that converges to a point $w \in \overline{q_f(U_i)} \subset W_f$. 
Note that by the above argument, we have another sequence
$\{w_n'\}$ in $q_f(U_i)$ that converges to $q_f(C)$.

Suppose that 
$w \neq q_f(C)$.
Then, there exist open neighborhoods $\Omega_0$ and $\Omega_1$
of $w$ and $q_f(C)$, respectively, in $W_f$ such that
$\Omega_0 \cap \Omega_1 = \emptyset$ and $\Omega_1 \subset q_f(U)$,
since $W_f$ is Hausdorff.
As $W_f$ is locally connected (see \cite[Theorem~8]{G5}),
we may assume that $\Omega_0$ and $\Omega_1$ are connected.
(Under these conditions, we can take $\Omega_0$ and $\Omega_1$
arbitrarily small.)
Note that there exist an integer $N$ such that for all $n \geq N$,
$w_n \in \Omega_0$ and $w'_n \in \Omega_1$.
Let us consider $\Omega_0 \cap q_f(U_i)$ and
$\Omega_1 \cap q_f(U_i)$ which contain $w_n$ and $w'_n$, respectively,
for all $n \geq N$. As both sequences converge to $0$
as sequences in $q_f(U_i)$ identified with $(0, 1)$, we see that both
$\Omega_0 \cap q_f(U_i)$ and
$\Omega_1 \cap q_f(U_i)$ have infinitely many connected
components, which are open in $W_f$.
Note that such connected components are contained
in $(0, 1)$ arbitrarily close to $0$, and hence each of them
can be considered to be an
open interval $(a, b)$ with $0 < a < b < 1$. As its closure in $W_f$
coincides with $[a, b]$, and $a$ and $b$ do not belong to $\Omega_0$ nor
$\Omega_1$ (since, otherwise, connected components would be strictly
larger), we see that the connected components of
$\Omega_0 \cap q_f(U_i)$ and
$\Omega_1 \cap q_f(U_i)$ are closed in $\Omega_0$
and in $\Omega_1$, respectively. 
This contradicts the connectedness of $\Omega_0$
or $\Omega_1$.
Note that such an argument is valid even if we replace $\Omega_0$ and $\Omega_1$
by smaller ones.
This is a contradiction.
This implies that $w = q_f(C)$ and hence
$r_i^{-1}$ is continuous.

Hence, $q_f(U_i \cup C) = q_f(U_i) \cup \{q_f(C)\}$ is 
homeomorphic to $[0, 1)$
in such a way that $q_f(C)$ corresponds to $0$.

These arguments show that 
we can naturally construct a continuous
surjection
$g : U \to S_\ell$ in such a way that
the image of each component of $U \setminus C$ coincides with
a component of $S_\ell \setminus \{p\}$
and that $g(C) = \{p\}$.

Then, we have a unique continuous map $\gamma : q_f(U) \to S_\ell$
such that $g = \gamma \circ q_f|_U$.
Note that $\gamma$ is a bijection.
Furthermore, its inverse is also continuous, since
$S_\ell$ is the union of $\ell$ copies of $[0, 1)$, say
$E_i$, $i = 1, 2, \ldots, \ell$, attached to $p$
along the points corresponding to $0$, each of $E_i$
is closed in $S_\ell$, and $\gamma^{-1}$ is continuous
on each $E_i$ as we saw above.
Therefore, $q_f(C)$ has the open neighborhood $q_f(U)$
homeomorphic to $S_\ell$.

Consequently, each point of $W_f$ has an open neighborhood
homeomorphic to $S_\ell$ for some $\ell \geq 0$.
(Note that when $\ell = 0$, the point is isolated,
and when $\ell = 2$, $S_2$ is homeomorphic to an open interval.)

Recall that $W_f$ is Hausdorff by \cite{G5}. 
By taking off, from $W_f$, the star-shaped small open
sets that are the $q_f$--images of
the finitely many star-like
open neighborhoods as above, we get
a compact topological $1$--dimensional manifold
with boundary. Thus, it is homeomorphic
to a finite disjoint union of copies of closed
intervals. As $W_f$ is the union of this $1$--dimensional
manifold and a finite number of star-shaped $1$--dimensional
complexes attached along the ``boundary points'',
we see that $W_f$ has the structure of a finite graph.

This completes the proof.
\end{proof}

\begin{rem}\label{rem0.7}
If $f : M \to \R$ is a smooth function
on a closed $1$--dimensional manifold, then
$M$ is a finite disjoint union of circles.
In this case, by an argument similar to
the above, we can show that the Reeb space
$W_f$ is a finite disjoint union of circles and points.
Thus, the Reeb space always has the structure of a graph.
\end{rem}

\begin{rem}\label{rem1}
Let $M$ be a closed manifold of dimension $m \geq 2$
and $f : M \to \R$ a smooth function with finitely
many critical values. Then, according to \cite{Sa22}
every critical contour
has a star-like open neighborhood.
Furthermore, there are at most finitely many critical contours.
Compare this with Example~\ref{ex2} of \S\ref{section3}, which
shows that there can be an infinite number of connected components of 
critical point set.
\end{rem}

\begin{rem}\label{rem2}
Let $M$ be a compact manifold of dimension $m \geq 2$ with
non-empty boundary and $f : M \to \R$ a smooth function.
For such a function, a \emph{critical contour}
is a contour which contains a critical point of
$f$ or $f|_{\partial M}$. If a contour is not critical,
then it is said to be \emph{regular}.
Under these terminologies, a \emph{cylindrical} open set and a 
\emph{star-like} open neighborhood of a critical contour can be
defined exactly the same way as in Definition~\ref{def2}.
Then, Theorem~\ref{thm2} holds for such smooth
functions on compact manifolds with non-empty boundary
as well. The proof is almost the same and is left to the reader.
\end{rem}

Theorem~\ref{thm2} can be generalized to
continuous functions on separable Peano continuum
as follows. Recall that a \emph{Peano continuum}
is a non-empty compact connected and locally connected
Hausdorff space.

\begin{dfn}\label{def3}
Let $f : X \to \R$ be a continuous function on a topological
space $X$. 

(1) A non-empty open, connected and saturated subset $U$ of $X$ is said to be
\emph{cylindrical} if
for every $z \in q_f(U)$,
$U \setminus q_f^{-1}(z)$ has exactly two connected components.

(2) An open, connected and saturated subset $U$ of $X$ is said to be
\emph{star-like} if for a contour $L \subset U$ of $f$, $U \setminus L$
consists of finitely many cylindrical subsets.
\end{dfn}

\begin{thm}\label{thm3}
Let $f : X \to \R$ be a continuous map
of a Peano continuum $X$ which is second countable. Then, the
Reeb space $W_f$ of $f$ has the structure of a finite graph
if and only if every contour has a star-like open neighborhood.
\end{thm}

Note that by Urysohn's metrization theorem,
the above $X$ is metrizable and separable.
Then, the above theorem can be proved by the same argument
as in the proof of Theorem~\ref{thm2} with
the help of some results obtained in \cite{G5}.
So, we omit the proof here.

\begin{rem}\label{rem4}
Instead of connected components of level
sets, if we consider path-connected components,
then the situation is very different.
See \cite[\S4]{Sa22}.
\end{rem}

\section{Examples}\label{section3}

In this section, we give various examples
of smooth functions on closed manifolds which have
interesting properties.

\begin{ex}\label{ex1}
Let $C \subset [0, 1]$ be an arbitrary closed subset
with $C \neq [0, 1]$.
By Whitney \cite{Wh1934}, there exists a smooth function $g : [0, 2]
\to [0, \infty)$ with $g^{-1}(0) = C \cup [1, 2]$.
Let $h : [0, 2] \to [0, \infty)$ be the smooth function defined by
$$h(t) = \int_0^t g(s)\, ds$$
for $t \in [0, 2]$.
We can easily show that the
critical point
set of $h$ coincides with $C \cup [1, 2]$.
(For the graph of $h$, refer to Fig.~\ref{fig13}.)

\begin{figure}[tbh]
\centering
\psfrag{2}{$2$}
\psfrag{1}{$1$}
\psfrag{0}{$0$}
\psfrag{h}{$h$}
\psfrag{t}{$t$}
\includegraphics[width=\linewidth,height=0.25\textheight,
keepaspectratio]{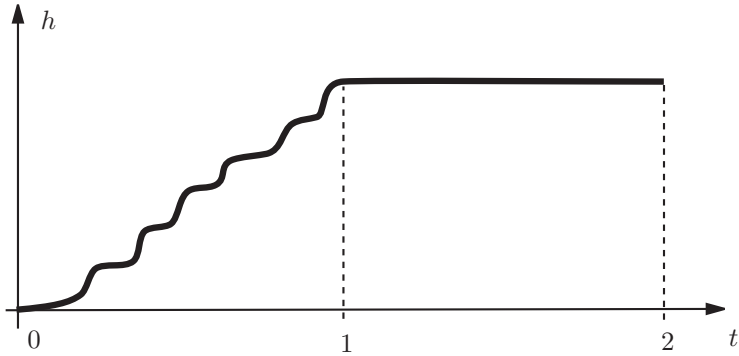}
\caption{Graph of the function $h$}
\label{fig13}
\end{figure}

For an appropriate $C$, $h$ has 
infinitely many critical values. For example,
consider the standard Cantor set in $[0, 1]$.

With the help of $h$, for every closed $m$--dimensional manifold 
$M$ with $m \geq 2$, we can construct a smooth function $f : M \to \R$
that has
infinitely many critical values
such that the Reeb space
$W_f$ has 
the structure of a finite graph.
For example, we first consider an arbitrary
Morse function $F$ on $M$, and around a maximal point $p$, we replace the
function with a smooth function of the form
$$H(x_1, x_2, \ldots, x_m) = h(2-(x_1^2 + x_2^2 + \cdots + x_m^2))+ c$$
with respect to certain local coordinates $(x_1, x_2, \ldots, x_m)$
around $p$ with $x_1^2 + x_2^2 + \cdots + x_m^2 < 2$
for some constant $c$.
Then, the Reeb space $W_f$ of the resulting function $f$ 
has the structure of a graph by Theorem~\ref{thm2}, since
every contour has a star-like open neighborhood.
In fact, $W_f$ is
homeomorphic to the Reeb space of $F$, which has the structure
of a finite graph.

Therefore, in Theorem~\ref{thm1} (\cite[Theorem~3.1]{Sa22}),
the condition that a smooth function should have
finitely many critical values is not a necessary condition for
the graph structure on the Reeb space.
\end{ex}

\begin{figure}[tbh]
\centering
\psfrag{1}{$1$}
\psfrag{0}{$0$}
\psfrag{-1}{$-1$}
\psfrag{f}{$f$}
\psfrag{S}{$S^2$}
\psfrag{R}{$\R$}
\includegraphics[width=\linewidth,height=0.6\textheight,
keepaspectratio]{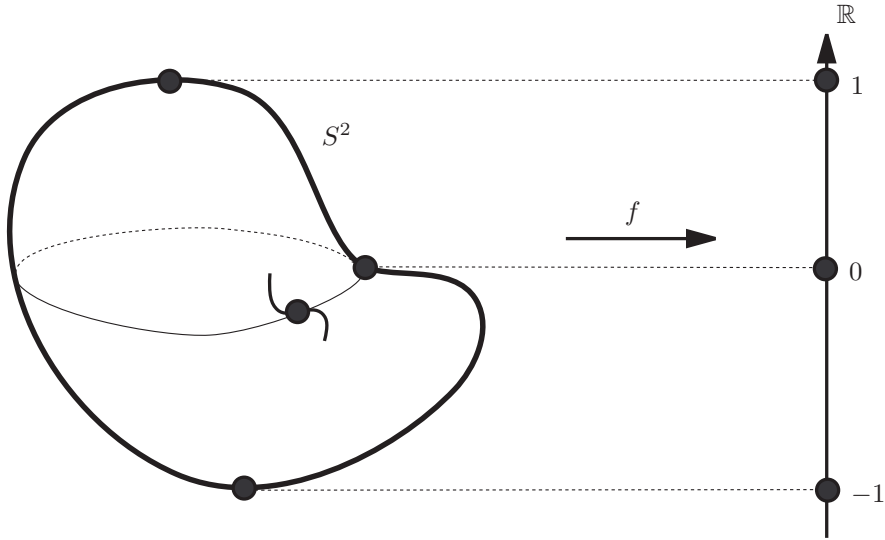}
\caption{A function $f : S^2 \to \R$ with infinitely many 
critical set components}
\label{fig14}
\end{figure}

\begin{ex}\label{ex2}
As we have seen in \cite{Sa22}, if a smooth function
$f$ on a closed manifold has only a finite number of
critical values, then the number of components
of a critical level set is always finite.

However, the critical point set itself may have infinitely
many components.
Let us give an example of a smooth
function on a closed manifold which has
infinitely many connected components of critical point set,
but which has only finitely many critical values.

We can construct a smooth function
$f : S^2 \to \R$ with the following properties, where
we identify $S^2$ with the unit sphere in $\R^3$.
\begin{enumerate}
\item The critical value set is $\{-1, 0, 1\}$.
\item $f^{-1}(-1) = \{(0, 0, -1)\}$ and
$f^{-1}(1) = \{(0, 0, 1)\}$.
\item $f^{-1}(0)$ coincides with the equator $S^2 \cap \{x_3=0\}$.
\item The critical point set in $f^{-1}(0) \cong S^1$
can be any prescribed closed subset $C$ of $S^1$.
In particular, it can be a Cantor set and
it may have uncountably many connected components.
\end{enumerate}

In fact, such an $f$ can be constructed in such
a way that $f$ is topologically equivalent to
the standard height function (see Fig.~\ref{fig14}).
More precisely, let $g : S^1 \to [0, 1]$ be a smooth
function such that $g^{-1}(0) = C$.
Then, we can construct the smooth function $f$
in such a way that the value of the 
derivative of $f$ with respect to
the upward unit tangent vector of $S^2$ at each point
on the equator
coincides with the $g$--value at the point. (At each point in $C$, we
arrange the function $f$ so that its second derivative
also vanishes, but the third derivative does not.)
\end{ex}

Compare the above example with \cite[Remark~0.3]{KP}.

\begin{ex}\label{ex3}
In \cite[Example~3.11]{Sa22}, we have constructed a smooth function
$f$ on each closed surface such that it has infinitely many
critical values, and that $W_f$ does not have the structure of a graph.
In the example, it was observed that each level set
has only finitely many connected components; however,
the number is not bounded.

In the following, we use the same notation as
in \cite[Example~3.11]{Sa22}.
Recall that
$D_n$ is the closed disk in $\R^2$ centered at the point 
$(1/n, 0$) with radius $1/2n(n + 1)$.
We warn the reader that
$$f_n(x) = g\left(2n(n + 1)(x-1/n)\right)$$
should be replaced by
$$f_n(x) = g\left(2n(n + 1)(x-(1/n, 0))\right)$$
for $x \in \R^2$
in \cite[Example~3.11]{Sa22}.

Set $d_n
= \sum_{k=n}^\infty c_k$ for $n \geq 1$.
Let us consider a 
smooth function $h : \R \to \R$
which is monotone increasing in the weak sense,
with the following properties.
\begin{enumerate}
\item $h(t) = 0$ for $t \leq 0$.
\item $h(t) = d_2$ for $t \geq 1$.
\item $h(t) = d_{n+1}$ for $t \in [1/n - 1/2n(n + 1),
1/n + 1/2n(n + 1)]$, for $n = 1, 2, 3, \ldots$.
\end{enumerate}

Then consider the smooth function
$\tilde{f} : \R^2 \to \R$ defined by
$$
\tilde{f}(x_1, x_2) = \begin{cases}
f(x_1, x_2) + h(x_1), & x_1 \leq 5/4, \\
f_1(x_1-1/2, x_2) + h(x_1), & x_1 \geq 5/4
\end{cases}
$$
for $x = (x_1, x_2) \in \R^2$.
For the graph of the function $\tilde{f}(x_1, 0)$, refer
to Fig.~\ref{fig23}.

\begin{figure}[tbh]
\centering
\psfrag{tf}{$\tilde{f}(x_1, 0)$}
\psfrag{x1}{$x_1$}
\psfrag{c1}{$d_1$}
\psfrag{c2}{$d_2$}
\psfrag{c3}{$d_3$}
\psfrag{1}{$1$}
\psfrag{a1}{$\frac{5}{4}$}
\psfrag{a2}{$\frac{3}{2}$}
\psfrag{a3}{$\frac{7}{4}$}
\psfrag{a4}{$\frac{3}{4}$}
\psfrag{a5}{$\frac{7}{12}$}
\psfrag{a6}{$\frac{1}{2}$}
\psfrag{a7}{$\frac{1}{3}$}
\includegraphics[width=\linewidth,height=0.4\textheight,
keepaspectratio]{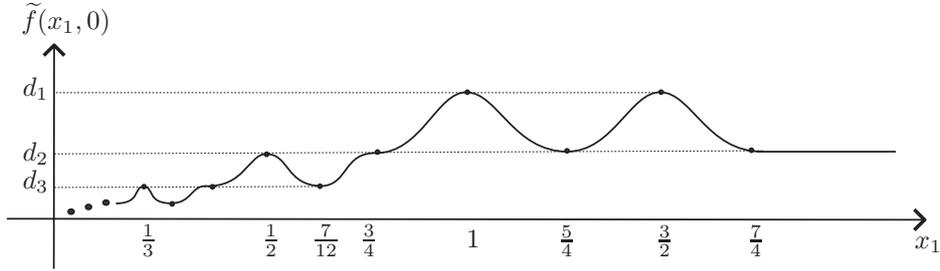}
\caption{Graph of the function $\tilde{f}(x_1, 0)$}
\label{fig23}
\end{figure}

With the help of $\tilde{f}$, we can construct a smooth function
$F : S^2 \to \R$ such that $F(S^2) =
[-d_1, d_1]$, having the following
properties:
\begin{enumerate}
\item $F^{-1}(t)$ has exactly two
components for $t \in [-d_1, d_1] - \{0\}$,
\item $F^{-1}(0)$ is connected, and
\item $W_F$ does not have the structure of a graph.
\end{enumerate}
See Figs.~\ref{fig20} and \ref{fig21}.
Note that the Reeb space $W_F$ does not have
the structure of a graph since the set of
``vertices'' does not have the discrete topology.
Note that $W_F$ is a so-called dendrite.
\end{ex}

\begin{figure}[tbh]
\centering
\psfrag{S}{$S^2$}
\psfrag{f}{$F$}
\psfrag{R}{$\R$}
\psfrag{a}{$d_1$}
\psfrag{b}{$-d_1$}
\psfrag{0}{$0$}
\psfrag{d2}{$d_2$}
\psfrag{md2}{$-d_2$}
\psfrag{d3}{$d_3$}
\psfrag{md3}{$-d_3$}
\includegraphics[width=\linewidth,height=0.5\textheight,
keepaspectratio]{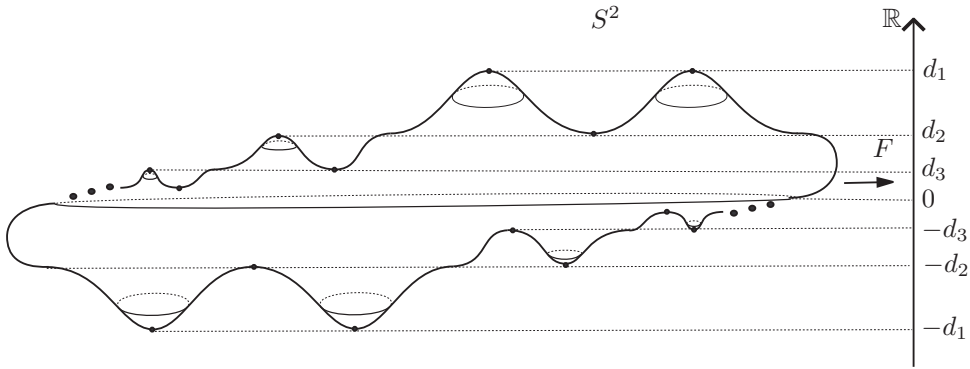}
\caption{A smooth function $F : S^2 \to \R$ with each level set
having exactly one or two connected components}
\label{fig20}
\end{figure}

\begin{figure}[tbh]
\centering
\psfrag{c1}{$d_1$}
\psfrag{0}{$0$}
\psfrag{c2}{$d_2$}
\psfrag{c3}{$d_3$}
\psfrag{c4}{$d_4$}
\psfrag{b}{$\overline{F}$}
\psfrag{W}{$W_F$}
\psfrag{R}{$\R$}
\psfrag{d1}{$-d_1$}
\psfrag{d2}{$-d_2$}
\psfrag{d3}{$-d_3$}
\psfrag{d4}{$-d_4$}
\includegraphics[width=\linewidth,height=0.4\textheight,
keepaspectratio]{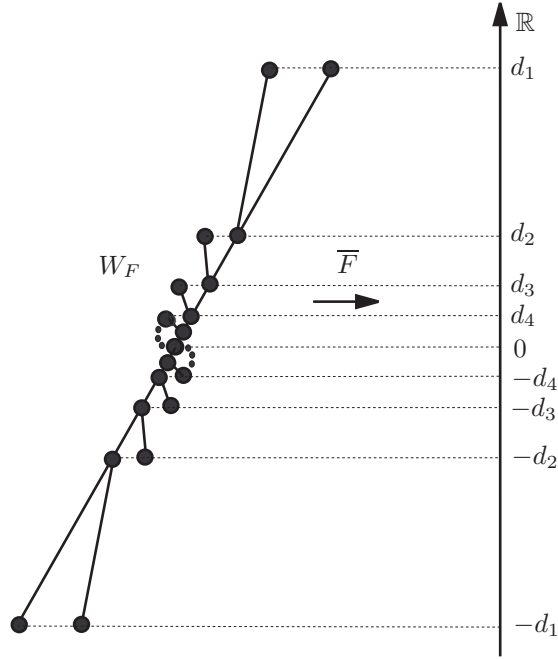}
\caption{Reeb space and Reeb function of $F$}
\label{fig21}
\end{figure}

\begin{ex}\label{ex4}
For $m \geq 2$, consider the standard Morse function $f : S^m \to \R$
defined by $f(x_1, x_2, \ldots, x_{m+1}) = x_{m+1}$,
where we identify $S^m$ with the unit sphere in $\R^{m+1}$.
Note that the Reeb space $W_f$ is homeomorphic to $[-1, 1]$
and $\bar{f} : W_f \to \R$ can be identified with the inclusion.
Let $\eta : [-1, 1] \to \R$ be any smooth function such that
each contour consists of a point. 
Then, the composition $g = \eta \circ q_f : S^m \to \R$
is a smooth function such that the Reeb space $W_g$
is identified with $W_f$ canonically homeomorphic to $[-1, 1]$,
and $\bar{g} : [-1, 1] \to \R$ is identified with $\eta$.
This means that even if the Reeb space has the structure
of a graph, the Reeb function restricted to each edge
may not be an embedding. Compare this with
Theorem~\ref{thm1}.
\end{ex}

Note also that in Theorem~\ref{thm1},
the graph structure of the
Reeb space 
satisfies that the vertices correspond exactly to
the connected components of level sets that contain critical points.
Therefore, the Reeb function restricted to each edge is an embedding
into $\R$ and
the resulting graph does not have a loop.
On the other hand, the graph structure in Theorem~\ref{thm2}
may have loops. For example, consider
the projection $S^1 \times S^1 \to S^1$ to the first factor
composed with a smooth function $S^1 \to \R$ such that
each contour consists of a point. Then, the resulting
smooth function $S^1 \times S^1 \to \R$ has a Reeb space
homeomorphic to $S^1$ which does have a structure of
a graph with a loop.

\section*{Acknowledgment}\label{ack}
The author would like to thank the organizers and
all the participants of the 17th International Workshop on
Real and Complex Singularities, held in S\~ao Carlos, Brazil,
in July 2022, which was an unforgettable event and
stimulated the author to write this paper.
The author has been supported in part by JSPS KAKENHI Grant Numbers 
JP22K18267, JP23H05437.

\end{document}